\documentclass[11pt]{amsart}
\usepackage{fullpage}
\usepackage{graphicx}
\usepackage{subfigure}
\usepackage{psfrag}
\usepackage{color}
\usepackage{amsmath}
\usepackage{amssymb}
\newtheorem{theorem}{Theorem}

\newtheorem{remark}[theorem]{Remark}

\title{A stabilized finite element formulation for advection--diffusion using the generalized finite element framework}

\author{D. Z. Turner} \author{K.~B.~Nakshatrala} \author{K.~D.~Hjelmstad}
\address{Correspondence to: Daniel Z. Turner, Department of Civil and Environmental Engineering, 
  2103 Newmark Laboratory, University of Illinois at Urbana-Champaign, Urbana, Illinois - 61801.}
\email{dzturne1@illinois.edu}
\address{Dr. Kalyana Babu Naskshatrala, Department of Civil and Environmental Engineering, 
  2524 Hydrosystems Laboratory, University of Illinois at Urbana-Champaign, Urbana, Illinois - 61801.}
\email{nakshatr@illinois.edu}
\address{Professor Keith D Hjelmstad, Department of Civil and Environmental Engineering, 
  3129e Newmark Laboratory, University of Illinois at Urbana-Champaign, Urbana, Illinois - 61801.}
\email{kdh@illinois.edu}

\date{\today}

\begin{document}


\begin{abstract}
	The following work presents a generalized (extended) finite element formulation for the advection--diffusion equation.  Using enrichment functions that represent the exponential nature of the exact solution, smooth numerical solutions are obtained for problems with steep gradients and high Peclet numbers (up to $Pe = 25$) in one and two-dimensions.  As opposed to traditional stabilized methods that require the construction of stability parameters and stabilization terms, the present work avoids numerical instabilities by improving the classical Galerkin solution with an enrichment function. To contextualize this method among other stabilized methods, we show by decomposition of the solution (in a multiscale manner) an equivalence to both Galerkin/least-squares type methods and those that use bubble functions.  This work also presents a strategy for constructing the enrichment function for problems with complex geometries by employing a global-local approach.
\end{abstract}


\maketitle

\section{INTRODUCTION}
%
The advection-diffusion equation represents transport by means of transmission and dissemination.  The problem is significant both in its application to physical phenomena and as a precursor to studying more complicated equations like the Navier--Stokes equations.  It is well known that the classical Galerkin method performs poorly for advection dominated transport problems. Spurious oscillations manifest in the solution due to the inherent truncation error of the discretized Galerkin approximation.  Numerous strategies have been proposed in the literature to overcome this deficiency, which include: those that add artificial diffusion or employ non-centered discretization of the advection operator (also called upwind schemes) \cite{SUPG, Heinrich, Brooks, SU, Johnson}, and multiscale models using bubble functions or wavelets \cite{Masud2, Dahmen}.  In many cases, an equivalence can be constructed between these methods \cite{Baiocchi, Brezzi3}.  For a complete review of stabilized methods for advection-diffusion, see \cite{Franca3} and the references therein.  

	Franca \emph{et al.} \cite{Franca3} show that for the advection-diffusion problem with small $\kappa$, small variations in the data produce large variations in the solution which represents the fundamental characteristic of instability.  The discrete classical Galerkin formulation also inherits the lack of stability of the continuous problem.  As the non-symmetric advection operator begins to dominate, the approximability of the Galerkin finite element method degrades.  For diffusion dominated flows, the solution is smooth, without sharp gradients or corners.  In such cases, the underlying finite element subspace adequately captures the exact solution.  Stabilized finite element methods improve the Galerkin formulation by reducing the truncation error (for example, Streamline-Upwind/Petrov-Galerkin) or increasing its approximability (for example, variational multiscale methods).  

  Generalized (or extended) finite element methods (GFEM or xFEM) have been successfully applied in a number of contexts, including: solid mechanics (cracks) \cite{Duarte, Simone, Sukumar}, fluid-structure interaction (moving interfaces) \cite{Gerstenberger}, immiscible multiphase flows \cite{Chessa}, and thermal problems with moving heat sources \cite{Merle}. Other mesh free methods include the method of finite spheres \cite{De} and the reproducing kernel method \cite{Liu}. This work presents the first application of the generalized finite element method to develop a stabilized-type finite element formulation.  To formulate a robust stabilized formulation requires the construction of stabilization parameters and stabilization terms.  In the generalized finite element context, one can improve the performance of the classical Galerkin formulation by designing an appropriate function that captures the solution better than the underlying finite element shape functions.  
  
  We show in a subsequent section that there are many similarities between generalized and stabilized finite element methods.  One can easily cast the generalized finite element method in such a way to obtain a stabilized formulation.  Additionally, the structure of the resulting stability parameter is identical to that of the Hughes variational multiscale method \cite{Masud2}.  The methods differ in that, whereas stabilized methods aim to add additional terms to the formulation, the underlying Galerkin formulation need not change for the generalized finite element method.  The methods also differ in that the generalized finite element method only enriches nodes in the vicinity of local phenomena such as a boundary layer.  Stabilized methods are typically applied to the entire domain.  Also, unlike stabilized methods, for generalized finite element methods the enforcement of Dirichlet boundary conditions is a bit more complicated as many of the enriched nodes fall on the boundary of the domain \cite{Cho}. Another complexity of generalized finite element methods regards numerical integration techniques for mesh free methods \cite{Chen}.

  The primary aim of this work is to produce smooth solutions to the advection--diffusion equation by means of the generalized finite element method rather than a stabilized formulation.  The numerical results presented show the successful application of the method to one and two-dimensional problems with high Peclet number (up to $Pe = 25$), with sharp corners, and boundary layer phenomena present. 
  
  In the sections that follow, we first introduce the governing equations and the classical Galerkin formulation.  We then discuss how the classical Galerkin formulation is enriched in the generalized finite element method.  Since the generalized finite element method is heavily dependent on the enrichment function, we then present several enrichment functions appropriate for the advection--diffusion equation.  Next, we discuss the enforcement of Dirichlet boundary conditions and the relationship between the generalized finite element method and stabilized methods.  Lastly, we present some numerical examples in one and two-dimensions.
  

  
%
  \section{GOVERNING EQUATIONS FOR ADVECTION--DIFFUSION}
 Let $\Omega$ be a bounded open domain, and $\Gamma$ ($\Gamma := \bar{\Omega} - \Omega$, where $\bar{\Omega}$ 
  is the closure of $\Omega$) be its piecewise smooth boundary. Let a scalar unknown field be denoted by 
  $u : \Omega \rightarrow \mathbb{R}$. Let the advection velocity field 
  be denoted by $\boldsymbol{\alpha}:\Omega \rightarrow \mathbb{R}^{nd}$, where $nd$ is the number of spatial dimensions and the diffusivity tensor be denoted as $\boldsymbol{\kappa}:\Omega \rightarrow \mathbb{R}^{nd \times nd}$. 
  As usual, $\Gamma$ is divided into two parts, denoted by 
  $\Gamma^{u}$ and $\Gamma^{t}$, such 
  that $\Gamma^{u} \cap \Gamma^{t} = \emptyset$ 
  and $\Gamma^{u} \cup \Gamma^{t}=\Gamma$.  $\Gamma^{u}$ is the part of the boundary 
  on which Dirichlet boundary conditions are prescribed, and $\Gamma^{t}$ 
  is the part of the boundary on which Neumann (or diffusive flux) boundary conditions are prescribed. 
  The governing equations for the advection-diffusion problem can be written as 
  \begin{align}
    \label{Eqn:AD_Equilibrium}
    \boldsymbol{\alpha} \cdot \nabla u - \nabla \cdot \boldsymbol{\kappa} \nabla u &= f \qquad \; \;  \ \mbox{in}  \quad \Omega \\
    \label{Eqn:AD_DirichletBC}
    u &= u^{\mathrm{p}} 
    \qquad \; \mbox{on} \quad \Gamma^{u} \\
    \label{Eqn:AD_NeumannBC}
    \boldsymbol{\kappa} \nabla u \cdot \boldsymbol{n} &= 
    t^{\mathrm{p}} \qquad \; \; \; \mbox{on} \quad \Gamma^{t} 
  \end{align}
  where $u^{\mathrm{p}}$ is the prescribed Dirichlet boundary condition, $t^{\mathrm{p}}$ is the prescribed flux, $\nabla$ is the gradient operator, $f$ is the volumetric source term, 
  and $\boldsymbol{n}$ is the unit outward normal vector 
  to $\Gamma$. The fluid is assumed incompressible, (i.e. $\nabla \cdot \boldsymbol{\alpha} = 0$). 
  
  The relative influence of advection and diffusion in equation \eqref{Eqn:AD_Equilibrium} is expressed by the nondimensional global Peclet number, $Pe = UL/\kappa$ where $L$ is the reference length, and $U$ is a reference velocity.  In the finite element context, we define the element Peclet number as $Pe^h = \parallel \boldsymbol{\alpha} \parallel h / 2 \kappa$, where $\parallel \cdot \parallel$ represents the $L_2$ norm and $h$ is the characteristic element length.
  \subsection{Classical Galerkin formulation}
  We now define function spaces that will be used in the remainder of the paper.  The functions spaces for the unknown scalar $u(\boldsymbol{x})$ and its associated weighting function $w(\boldsymbol{x})$ are 
\begin{align}
  \label{Eqn:AD_Function_Space_for_u}
  \mathcal{U} &:= \{u  \; \big| \; u \in H^{1}(\Omega),
  u = u^{\mathrm{p}} \; \mathrm{on} \; \Gamma^u \} \\
  \label{Eqn:AD_Function_Space_for_w}
  \mathcal{W} &:= \{w  \; \big| \; w \in H^{1}(\Omega),
  w = 0 \; \mathrm{on} \; \Gamma^u \}
\end{align}
The classical Galerkin formulation for the advection-diffusion problem is written as: Find $u(\boldsymbol{x}) \in \mathcal{U}$ such that 
\begin{align}
  \label{Eqn:AD_Classical_Functional}
  a(w;u) + 
  b(w;u) = 
  l(w) & \quad \forall \ w \in \mathcal{W}
\end{align}
Let us define the bilinear forms as: 
\begin{align}
  \label{Eqn:AD_a_Galerkin}
  a(w;u) &:= 
  \int_{\Omega} w(\boldsymbol{\alpha} \cdot \nabla u) \; \mathrm{d} \Omega \\
  \label{Eqn:AD_b_Galerkin}
  b(w;u) &:= 
  \int_{\Omega}  \nabla w \cdot \boldsymbol{\kappa} \nabla u \; \mathrm{d} \Omega  
\end{align}
and the linear functional as 
\begin{align}
  \label{Eqn:AD_f_Galerkin}
  l(w) & := \int_{\Omega} wf \; \mathrm{d} \Omega + \int_{\Gamma^t} wt \; \mathrm{d} \Gamma
\end{align}
Once the weak formulation of the governing equations is established, the approximate solution based on the finite element method is determined in the usual manner.  First one chooses the approximating finite element spaces, which (for a conforming formulation) will be the finite element function spaces for the unknown scalar $u$ and the weighting function associated with $u$ denoted by $\mathcal{U}^{h} 
\subseteq \mathcal{U}$ and $\mathcal{W}^h \subseteq \mathcal{W}$ respectively.  The finite element formulation of the weak form reads: Find $u^h(\boldsymbol{x}) \in \mathcal{U}^h$ such that
\begin{align}
  \label{Eqn:AD_FEM}
  a(w^h;u^h) + 
  b(w^h;u^h) = 
  l(w^h) & \quad \forall \ w^h \in \mathcal{W}^h 
\end{align}
\subsection{Stabilized methods}
It is well known that the standard Galerkin method performs poorly for advection-dominated flow (i.e. $Pe^h  \gg 1$) \cite{Gresho}. Excessive mesh refinement is necessary to avoid spurious oscillations that propagate throughout the solution.  One popular method for dealing with such oscillations involves the use of \emph{stabilized methods} (For examples of stabilized finite element methods for the advection-diffusion equation see \cite{Gresho,Mulder,Donea,SU,SUPG,Masud2,Brezzi3} and the references therein).  \emph{Stabilized methods} perturb the Galerkin formulation with mesh dependent terms.  In some cases, the perturbation is designed in an ad hoc manner, while for other methods the perturbation preserves consistency by being residual based.  In the present context, we introduce a stabilized-type formulation that is based on improving the classical Galerkin approximation by enrichment with a function that captures certain features of the solution.

%
  \section{ENRICHED APPROXIMATION}
In this section we present a generalized finite element approach to solving the advection-diffusion equation by employing the partition of unity framework as detailed in \cite{Melenk, Simone, Oden, Taylor}.
Let us divide the domain $\Omega^h$ into ``$Nele$'' non-overlapping subdomains 
$\Omega^h_e$ (which in the finite element context will be elements) such 
that
%
\begin{equation}
  \Omega^h = \overset{Nele}{\underset{e = 1}{\bigcup}} \Omega^h_{e}
\end{equation}
The boundary of element $\Omega^h_{e}$ is denoted by $\Gamma^h_{e}$. We then establish the standard finite element basis, $\mathcal{U}^h_0$ as the space of complete polynomials $P^k(\Omega^h_{e})$ of order $\leq k$ over each element
%
\begin{align}
  \label{Eqn:AD_function_space_for_Uh}
  \mathcal{U}^h_0 := \mathrm{span}\{N_i\}^{F}_{i=1} \; \mathrm{where} \;  \{ N_i \in \left[ C^0(\Omega^h_{e}) \right]^{nd} \ : \ N_i\bigr| \Omega^h_{e} \in \left[ P^k(\Omega^h_{e}) \right]^{nd} \; \mathrm{and} \;
  N_i\bigr|_{\Gamma^u} = 0\} 
\end{align}
  where $N_i$ is often referred to as the nodal shape functions, and $F$ is the number of nodal finite element shape functions. In the numerical examples, piecewise linear approximations are used for $N$.
  Consider a set of overlapping arbitrarily shaped domains $\omega_{i}$ centered at node $i$ that define the support for each shape function, $N_i$.  Also consider enrichment functions $\mathcal{H}_m$ that represent special approximating properties in the proximity of local phenomena.  A generalized finite element basis for $\mathcal{U}^h_0$ that preserves a conforming approximation can be constructed as
\begin{align}
  \label{Eqn:AD_function_space_for_Uh_2}
  \mathcal{U}^h_0 := \mathrm{span} \left\{ \{ N_i \} \cup \{ N_i\mathcal{H}_m \}^{E}_{m=1} \right\}^{F}_{i=1} 
\end{align}
where $E$ represents the number of enrichment functions for node $i$.  The unknown scalar field $u(\boldsymbol{x})$ may then be approximated as
\begin{align}
  \label{Eqn:AD_App_u}
  u(\boldsymbol{x}) = \sum_i \bar{u}_i N_i(\boldsymbol{x}) + \sum_i N_i(\boldsymbol{x})\left( \sum_m (u'_i)_m \mathcal{H}_m(\boldsymbol{x})  \right)
\end{align}
where $\bar{u}_i$ represents the extra coefficients on the enriched nodes.
\begin{remark}
Typically, when recovering the solution for a standard finite element analysis, the coefficients $u_i$ are equal to the value of the unknown scalar field $u(\boldsymbol{x})$ at the nodes.  In the case of generalized finite elements, the coefficients $u_i$ are not equal to the unknown scalar field at the nodes.
\end{remark}
Since our focus is on capturing the phenomena local to the boundary of the domain of interest, we need not enrich all the nodes in the domain.  For the case in which the finite element space is enriched by only a single function $\mathcal{H}(\boldsymbol{x})$ at nodes in the local support of the boundary, the unknown scalar field $u(\boldsymbol{x})$ may be approximated by
\begin{align}
  \label{Eqn:AD_App_u_single}
  u(\boldsymbol{x}) = \sum_{i \in I} \bar{u}_i N_i(\boldsymbol{x}) + \sum_{j \in J} u'_j N_j(\boldsymbol{x})  \mathcal{H}(\boldsymbol{x})
\end{align}
where $I$ represents all nodes in the mesh and $J$ represents the nodes that form a partition of unity for the function $\mathcal{H}(\boldsymbol{x})$
\begin{align}
  \label{Eqn:AD_define_J}
  J := \{ j \in I \ : \ \omega_j \cap \Gamma_I \neq \emptyset \}
\end{align}
To simplify the proceeding steps, we introduce the vectors of unknowns at each node, $\boldsymbol{u}$ and element shape functions $\boldsymbol{N}$,
\begin{align}
  \label{Eqn:AD_define_hats}
  \boldsymbol{u}^T = [ \bar{u}_1, \; \cdots, \; \bar{u}_n, \; u'_1, \; \cdots, \; u'_n] \; ; \; \boldsymbol{N}(\boldsymbol{x}) = [ N_1, \; \cdots, \; N_n, \; \mathcal{H}N_1, \; \cdots, \; \mathcal{H}N_n]
\end{align}
where $n$ is the number of nodes per element, such that
\begin{align}
\label{Eqn:AD_combined_shapes}
	u(\boldsymbol{x}) = \boldsymbol{N}(x) \boldsymbol{u}
\end{align}
Note that for elements that do not have any enriched nodes, $\boldsymbol{N} = [N_1, \; \cdots, \; N_n]$ and $\boldsymbol{u}^T = [\bar{u}_1, \; \cdots, \; \bar{u}_n]$.
%
After substitution of equation \eqref{Eqn:AD_combined_shapes} into equation \eqref{Eqn:AD_FEM} and noting the arbitrariness of $w$, the discrete form of the generalized finite element method of the advection-diffusion equations may be written as the assembled sum contribution of each element,
\begin{align}
  \label{Eqn:AD_Discrete_Functional}
  \mathbb{A} \sum_{e = 1}^{Nele} a(\boldsymbol{N};\boldsymbol{N})\boldsymbol{u} + 
  b(\boldsymbol{N};\boldsymbol{N})\boldsymbol{u} = 
  \mathbb{A} \sum_{e = 1}^{Nele} f(\boldsymbol{N})
\end{align}
where $\mathbb{A}$ is the standard assembly operator.  Equation \eqref{Eqn:AD_Discrete_Functional} is written in the familiar matrix form as
\begin{align}
  \label{Eqn:AD_Matrix_Functional}
  \boldsymbol{K} \boldsymbol{u} = \boldsymbol{f}
\end{align}
where
\begin{align}
  \label{Eqn:AD_Matrix_Functional_Define}
  \boldsymbol{K} &= \mathbb{A} \sum_{e = 1}^{Nele} a(\boldsymbol{N};\boldsymbol{N}) + 
  b(\boldsymbol{N};\boldsymbol{N}) \\
  \boldsymbol{f} &= \mathbb{A} \sum_{e = 1}^{Nele} f(\boldsymbol{N})
\end{align}
Equation \eqref{Eqn:AD_Matrix_Functional} represents a linear system of equations with the number of unknowns equal to the number of nodes plus the number of enriched nodes.
%
  \section{DESIGN OF ENRICHMENT FOR ADVECTION--DIFFUSION}
  Typically, the enrichment function or set of functions $\mathcal{H}(\boldsymbol{x})$ are chosen to provide an improved estimate for the solution in the vicinity of local phenomena.  In some cases, the exact solution is known in the immediate vicinity or at least the nature of the exact solution is known.  In such cases, $\mathcal{H}(\boldsymbol{x})$ may be choses as a particular part of the solution.  For complex geometries, or for problems for which little is known about the solution, we introduce a global-local approach to choosing the enrichment function. In a subsequent section, we explore the dependence of the solution on the enrichment function selected.
  \subsection{One-dimensional problem}
\begin{figure}[htb!]
	\centering
  \includegraphics[scale=0.35]{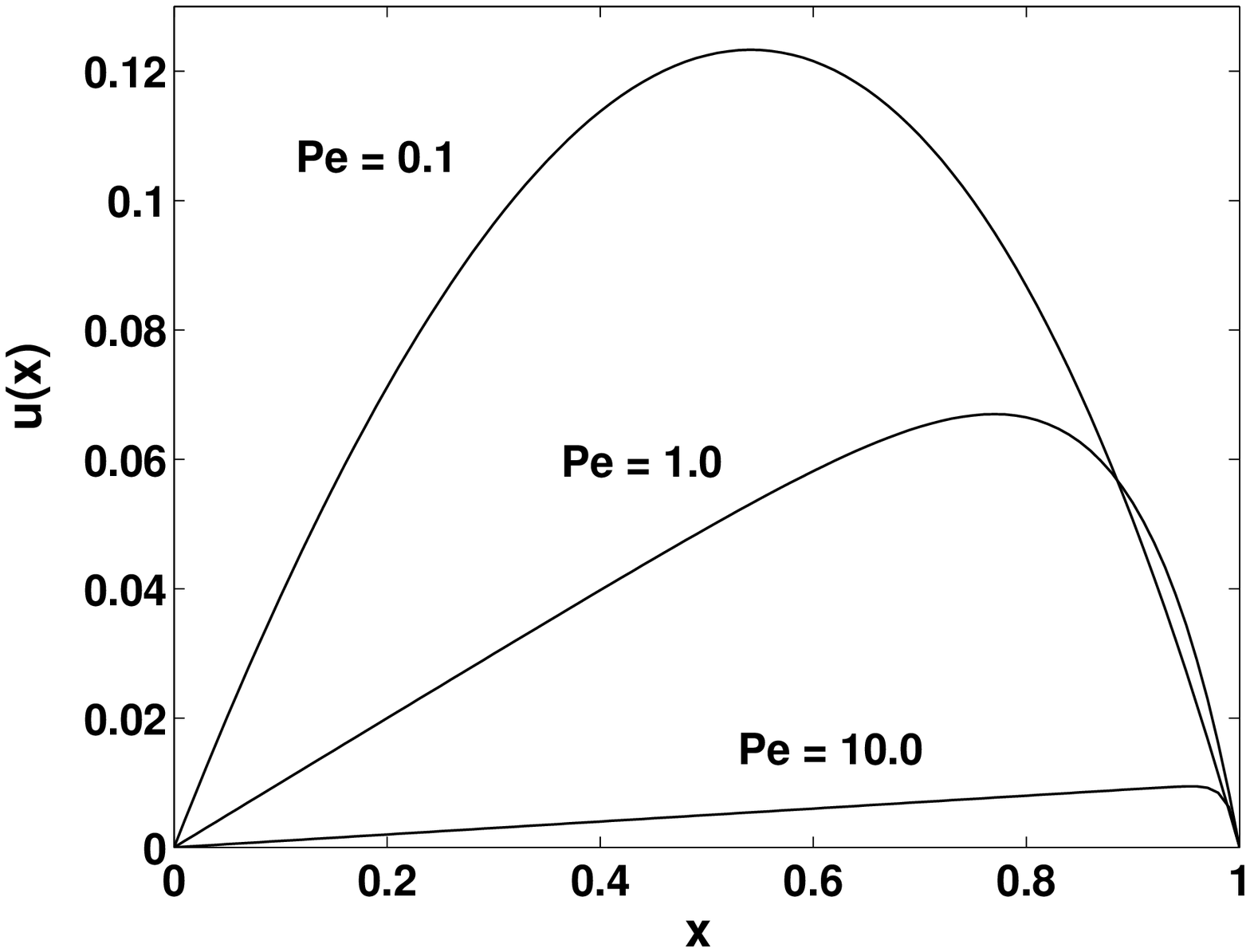}
	\caption{One-dimensional advection-diffusion: exact solution for various $Pe$.}
	\label{fig:1DExactSol}
\end{figure}

  Consider the advection-diffusion problem in one dimension on a domain of unit length with the following boundary conditions
  \begin{align}
    \label{Eqn:AD_1D_ad_dif}
    \alpha  \frac{\partial u}{\partial x} - \kappa \frac{\partial^2u}{\partial x^2} &= 1 \qquad \; \;  \ \mbox{in}  \quad \Omega \\
    \label{Eqn:AD_ad_dif_DirichletBC}
    u(x) &= 0  \qquad \; \;  \ \mbox{on}  \quad \Gamma
  \end{align}
Assuming a constant advection speed, $\alpha$, and diffusivity, $\kappa > 0$, the exact solution to this problem is
  \begin{align}
    \label{Eqn:AD_1D_exact}
    u(x) = \frac{1}{\alpha} \left( x - \frac{1 - e^{\alpha x / \kappa}}{1 - e^{\alpha / \kappa}} \right)
  \end{align}
The nature of the solution is highly dependent on the Peclet number, $Pe^h$.  For high $Pe^h$ a thin boundary layer is present toward the outflow boundary.  The exact solution for increasing $Pe^h$ is shown in Figure \ref{fig:1DExactSol}.  The simplest enrichment function, motivated by the exact solution, is to use
\begin{align}
 \label{Eqn:AD_H_a}
 \mathcal{H}_a(x) = e^{\gamma x}
\end{align}
where $\gamma = \alpha / \kappa$. In this case, we take advantage of the exponential nature of the solution.  Analogously, in linear fracture mechanics problems, one knows that the stress field grows as $1/\sqrt{r}$, where $r$ is the radius from the crack tip.  For advection-diffusion we know that the solution has an exponential nature.

\begin{figure}[htb!]
	\centering
  \includegraphics[scale=0.35]{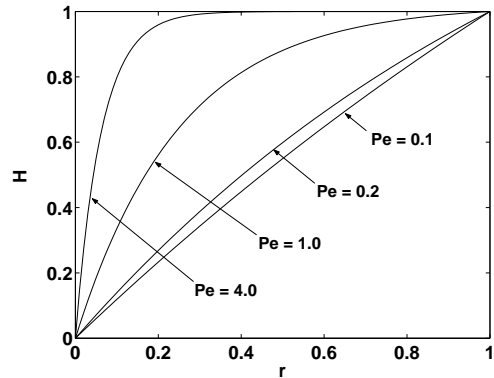}
	\caption{1D enrichment function as a function of the radius $r = 1 - x$ from the outflow boundary.}
	\label{fig:1DH}
\end{figure}

Of primary interest for flow problems is the application of the generalized finite method to capturing the localized nature in the boundary layer.  As such, the enforcement of Dirichlet boundary conditions becomes more involved as many of the enriched nodes lie on the Dirichlet boundary.  The enforcement of Dirichlet boundary conditions is discussed in detail in the next section, but we wish to construct here an alternative enrichment function that does not interfere with the Dirichlet boundary.  Building on the enrichment function $\mathcal{H}_a(x)$ as motivated by the exact solution, we can construct an alternate normalized enrichment function $\mathcal{H}_b(x)$ that vanishes on the boundary,
\begin{align}
 \label{Eqn:AD_H_b}
 \mathcal{H}_b(x) = 1 - \frac{1 - e^{\gamma x}}{1 - e^{\gamma}}
\end{align}
where again $\gamma = \alpha / \kappa $.  Figure \ref{fig:1DH} shows $\mathcal{H}_b(x)$ for various $Pe$.

Closer examination of equation \eqref{Eqn:AD_H_b} reveals that $\mathcal{H}_b(x)$ is not well-defined numerically for $\gamma = 0$ or $\gamma$ much greater than $800$, which represent the cases of pure diffusion or advection dominated flows respectively.  An alternative enrichment $\mathcal{H}_c(x)$ may be selected that is a polynomial representation of the exponential nature of $\mathcal{H}_b(x)$
\begin{align}
 \label{Eqn:AD_H_c}
 \mathcal{H}_c(x) = 1 - x^{\gamma}
\end{align}
$\mathcal{H}_c(x)$ both vanishes at the outflow boundary and is well-defined for pure diffusion and advection dominated flows.
 \subsection{Global-local enrichment function}
 In multi-dimensional problems on complex geometries, choosing an appropriate enrichment function becomes more difficult.  As such, one can construct an enrichment function based on the solution to the problem at low $Pe^h$, for which the classical Galerkin formulation is stable.  In the global-local approach,  one begins with a stable solution at low $Pe^h$ and iteratively increases the $Pe^h$, using the previous solution as the enrichment function.  The algorithm goes as follows,\\

\begin{tabular}{l}
\hline
Global-local solution strategy \\
\hline 
Solve the Galerkin problem (equation \eqref{Eqn:AD_FEM}) for $u^h_0$, with $Pe^h_0 = 1.0$ \\
$\mathcal{H}_0 = u^h_0$ \\
$\Delta Pe^h = (Pe^h_{N} - Pe^h_{0}) / N$, where $N$ is the number of iterations\\
for $i=1$ to $N$ \\
$\qquad \mathcal{H}_i = u^h_{i-1}$, $Pe^h_i = Pe^h_{i-1} + \Delta Pe^h$\\
$\qquad$ Solve the GFEM problem (equation \eqref{Eqn:AD_Matrix_Functional}) for $u^h_i$, with $Pe = Pe^h_i$\\
end\\
\hline
 \end{tabular} \\
 \\
 The preceding solution strategy may also be considered as a continuation method in the context of flow problems \cite{Deuflhard}.  The aim of a continuation method is to solve the problem for a high characteristic flow speed using information from the solution at a lower characteristic speed.
%
\section{ENFORCEMENT OF DIRICHLET BOUNDARY CONDITIONS}
Among the various methods available in mesh-free methods for the enforcement of Dirichlet boundary conditions, one of the simplest ways to implement is the penalty method.  The penalty method introduces no additional unknowns and preserves the banded structure, positive definiteness, and symmetry of the stiffness matrix.
\subsection{Penalty enforcement}  
  We begin with a variational statement from which the advection-diffusion equation may easily be derived in one dimension.  Consider the variational functional $\Pi(u)$
\begin{align}
    \label{Eqn:AD_Variational_Form}
	\Pi(u) = \int_{\Omega} \frac{\kappa}{2} \parallel \nabla u \parallel^2 e^{-\alpha x / \kappa} - u(x) f e^{-\alpha x / \kappa} \; \mathrm{d} \Omega 
\end{align}
  where $\nabla$ is the gradient operator, $\parallel \cdot \parallel$ represents the appropriate norm, and 
 $f$ represents the body force.  Solving the advection-diffusion equation is equivalent to extremizing $\Pi(u)$ subject to $u(x) = u_0$ on $\Gamma^u$.  Applying the stationarity condition to $\Pi(u)$,
\begin{align}
    \label{Eqn:AD_D_I}
	\delta \Pi = \frac{\partial \Pi}{\partial u(x)} - \frac{\mathrm{d}}{\mathrm{d}x}\left(\frac{\partial\Pi}{\partial \nabla u}\right) = 0
\end{align}
we have
\begin{align}
    \label{Eqn:AD_D_I2}
	\int_{\Omega} \alpha \frac{\partial u}{\partial x} - \kappa \frac{\partial^2 u}{\partial x^2} - f  \; \mathrm{d} \Omega = 0 \; ; \; \mathrm{s.t.} \; u(x) = u_0 \; \mathrm{on} \; \Gamma^u
\end{align}
Using a penalty method, we change the constrained extremization problem to an unconstrained extremization problem by adding a penalty functional, $\lambda \Lambda(u) / 2$ to $\Pi(u)$
\begin{align}
    \label{Eqn:AD_Penalty_Form}
	\Lambda(u) = \int_{\Gamma^u} \frac{\lambda}{2} (u(x) - u_0)^2 \; \mathrm{d} \Gamma 
\end{align}
where $\lambda$ is the penalty parameter.  Extermizing the sum of both functionals $\Pi(u)$ and $\Lambda(u)$, the solution will satisfy the Dirichlet boundary conditions as $\lambda$ approaches $\infty$.  The resulting functional, that incorporates a weak enforcement of the essential boundary conditions is written:
\begin{align}
    \label{Eqn:AD_Total_Penalty}
	\int_{\Omega} \alpha \frac{\partial u}{\partial x} - \kappa \frac{\partial^2 u}{\partial x^2} - f  \; \mathrm{d} \Omega + \lambda \int_{\Gamma^u} (u(x) - u_0) \; \mathrm{d} \Gamma = 0 
\end{align}
For mesh-free methods, the enforcement of essential boundary conditions must be handled with great care to avoid over-constraint phenomena.  To eliminate spurious results when enforcing the Dirichlet boundary conditions weakly, one can use trapezoidal numerical integration rather than Gauss integration to evaluate the boundary terms ($\int_{\Gamma^u} (\cdot) \; \mathrm{d}\Gamma$) in equation \eqref{Eqn:AD_Total_Penalty} \cite{Cho}.  
When enrichments are used that do not vanish on the Dirichlet boundary (for example $\mathcal{H}_a$), the boundary conditions will be enforced in a weak fashion according to equation \eqref{Eqn:AD_Total_Penalty}.
%
  \section{GFEM AND STABILIZED METHODS} 
 Due to the proliferation of stabilized methods for the advection--diffusion problem, as pointed out in the introduction, we wish to cast the generalized finite element formulation for advection--diffusion in a different manner to explore its relationship with stabilized methods. To explore this relationship, we begin by establishing a variational mulitiscale framework, in which we separate the problem into two sub-problems: the first representing the standard finite element solution, and the second the enriched finite element solution.  We then solve the enriched problem in terms of the standard variables and substitute back into the standard finite element problem.  The result is a stabilized formulation equivalent to Galerkin/least-squares (GLS) or streamline-upwind/Petrov-Galerkin (SUPG) and the structure of the stabilization parameter is identical to that of the Hughes variational multiscale method.  An equivalence between stabilized methods such as GLS and multiscale methods that employ bubble functions has previously been developed in works such as \cite{Baiocchi, Brezzi3}, but we wish to extend this equivalence to the generalized finite element context.
 
 Consider the advection--diffusion equation enriched with a single function, $\mathcal{H}(\boldsymbol{x})$ and constant diffusivity $\kappa$. We first decompose the unknown field in equation \eqref{Eqn:AD_App_u_single} and its weighting function into a standard Galerkin part, $\bar{\boldsymbol{u}}$, and an enriched part, $\boldsymbol{u}'$, as follows:
 \begin{align}
 \label{Eqn:Decompose_Solution}
 	u(\boldsymbol{x}) = \bar{u}(\boldsymbol{x}) + u'(\boldsymbol{x}) \quad  ;  \quad  w(\boldsymbol{x}) = \bar{w}(\boldsymbol{x}) + w'(\boldsymbol{x})
 \end{align}
 where
 \begin{align}
 \label{Eqn:Multiscale_Defs_1}
 	\bar{u}(\boldsymbol{x}) = \sum_{i \in I} \bar{u}_i N_i(\boldsymbol{x}) = \boldsymbol{N}\bar{\boldsymbol{u}} \quad  &; \quad u'(\boldsymbol{x}) = \sum_{j \in J} u'_j N_j(\boldsymbol{x}) \mathcal{H}(\boldsymbol{x})  =  \boldsymbol{N}'\boldsymbol{u}'\\
	\label{Eqn:Multiscale_Defs_2}
 	\bar{w}(\boldsymbol{x}) = \sum_{i \in I} \bar{w}_i N_i(\boldsymbol{x}) = \boldsymbol{N}\bar{\boldsymbol{w}}  \quad &; \quad 
w'(\boldsymbol{x}) = \sum_{j \in J} w'_j N_j(\boldsymbol{x}) \mathcal{H}(\boldsymbol{x}) = \boldsymbol{N}'\boldsymbol{w}'
 \end{align}
 The weak form (equation \eqref{Eqn:AD_FEM}) may now be written as
 \begin{align}
 \label{Eqn:Split_Weak_Total}
 	a(\bar{w} + w';\bar{u} + u') + 
  	b(\bar{w} + w';\bar{u} + u') = 
  	l(\bar{w} + w') 
 \end{align} 
 Due to the linearity of the weighting function and the trial solution, we can decompose \eqref{Eqn:Split_Weak_Total} into two subproblems: the \emph{Galerkin subproblem},
 \begin{align}
  \label{Eqn:Galerkin_Subproblem}
 	a(\bar{w};\bar{u}) + a(\bar{w};u') +
  	b(\bar{w};\bar{u}) + b(\bar{w};u') = 
  	l(\bar{w}) 
 \end{align}
 and the \emph{enriched subproblem}
 \begin{align}
  \label{Eqn:Enriched_Subproblem}
 	a(w';\bar{u}) + a(w';u') +
  	b(w';\bar{u}) + b(w';u') = 
  	l(w') 
 \end{align}
 Using equation \eqref{Eqn:AD_Function_Space_for_w} and integration by parts as follows,
 \begin{align}
 	\int_{\Omega} w(\boldsymbol{\alpha} \cdot \nabla u) \; \mathrm{d} \Omega & = - \int_{\Omega} \boldsymbol{\alpha} \cdot \nabla w \cdot u \; \mathrm{d} \Omega = c(w;u) \\
	 \int_{\Omega}  w \nabla \cdot \boldsymbol{\kappa} \nabla u \; \mathrm{d} \Omega & = -\int_{\Omega} w \boldsymbol{\kappa}\nabla^2 u  \; \mathrm{d} \Omega = d(w;u)
 \end{align}
 The \emph{Galerkin subproblem} and \emph{enriched subproblem} may alternatively be written as
 \begin{align}
  \label{Eqn:Galerkin_Subproblem_Alt}
 	a(\bar{w};\bar{u}) + 
  	b(\bar{w};\bar{u}) + c(\bar{w};u') + d(\bar{w};u') = 
  	l(\bar{w}) 
 \end{align}
and
 \begin{align}
  \label{Eqn:Enriched_Subproblem_Alt}
 	a(w';\bar{u}) + a(w';u') +
  	b(w';u') + d(w';\bar{u}) = 
  	l(w') 
 \end{align} 
respectively.

Consider the \emph{enriched subproblem}, equation \eqref{Eqn:Enriched_Subproblem_Alt}, which can be written as
\begin{align}
	\int_{\Omega} w'(\boldsymbol{\alpha} \cdot \nabla u') + \nabla w' \cdot \kappa \nabla u' \; \mathrm{d} \Omega = - \int_{\Omega} w' \bar{r} \; \mathrm{d}\Omega
\end{align}
where $\bar{r} = -\boldsymbol{\alpha} \cdot \nabla \bar{u}  + \kappa \nabla^2 \bar{u} - f$.  Using equations \eqref{Eqn:Multiscale_Defs_1} and \eqref{Eqn:Multiscale_Defs_2}, the above equation can be expressed in discrete form as
\begin{align}
	\boldsymbol{w}'^T \int_{\Omega} \boldsymbol{N}' \boldsymbol{\alpha}^T \nabla \boldsymbol{N}' + \kappa \nabla\boldsymbol{N}'^T\nabla\boldsymbol{N}' \; \mathrm{d}\Omega \boldsymbol{u}' = - \boldsymbol{w}'^T \int_{\Omega} \boldsymbol{N}'^T \bar{r} \; \mathrm{d}\Omega
\end{align}
Noting the arbitrariness of $\boldsymbol{w}'$, we have
\begin{align}
	\boldsymbol{u}' = - \frac{\int_{\Omega} \boldsymbol{N}'^T \bar{r} \;  \mathrm{d} \Omega}{ \int_{\Omega} \boldsymbol{N}'^T \boldsymbol{\alpha}^T \nabla \boldsymbol{N}' + \kappa \nabla \boldsymbol{N}'^T \nabla \boldsymbol{N}' \; \mathrm{d} \Omega }
\end{align}
Due to equation \eqref{Eqn:Multiscale_Defs_1} the enriched solution can be written as a function of the variables of the \emph{Galerkin subproblem} as
\begin{align}
	u'(\boldsymbol{x}) = - \frac{\boldsymbol{N}' \int_{\Omega} \boldsymbol{N}'^T \;  \mathrm{d} \Omega  \bar{r}}{ \int_{\Omega} \boldsymbol{N}'^T \boldsymbol{\alpha}^T \nabla \boldsymbol{N}' + \kappa \nabla \boldsymbol{N}'^T \nabla \boldsymbol{N}' \; \mathrm{d} \Omega } = -\tau \bar{r}
\end{align}
where we have made the assumption that $\bar{r}$ is constant over the element domain, which is exactly true for a constant body force using linear elements, and we have introduced the familiar intrinsic time scale, $\tau$, also known as the stabilization parameter
\begin{align}
\label{Eqn:Define_Tau}
\tau = \frac{\boldsymbol{N}' \int_{\Omega} \boldsymbol{N}'^T \;  \mathrm{d} \Omega}{ \int_{\Omega} \boldsymbol{N}'^T \boldsymbol{\alpha}^T \nabla \boldsymbol{N}' + \kappa \nabla \boldsymbol{N}'^T \nabla \boldsymbol{N}' \; \mathrm{d} \Omega }
\end{align}
Substituting equation \eqref{Eqn:Define_Tau} into the \emph{Galerkin subproblem}, equation \eqref{Eqn:Galerkin_Subproblem_Alt}, we have
 \begin{align}
  \label{Eqn:Galerkin_Subproblem_Stab}
 	a(\bar{w};\bar{u}) + 
  	b(\bar{w};\bar{u}) + c(\bar{w};\tau \bar{r}) + d(\bar{w};\tau \bar{r}) = 
  	l(\bar{w}) 
 \end{align}
Using the general inner product, $(w;u) = \int_{\Omega} (\cdot) \mathrm{d} \Omega$ Equation \eqref{Eqn:Galerkin_Subproblem_Stab} can be written more explicitly as
\begin{align}
	(\bar{w}, \boldsymbol{\alpha} \cdot \nabla \bar{u}) + (\nabla \bar{w}, \kappa \nabla \bar{u}) + (\boldsymbol{\alpha} \cdot \nabla \bar{w} + \kappa \nabla^2 \bar{w}, \tau(\boldsymbol{\alpha} \cdot \nabla \bar{u} - \kappa \nabla^2 \bar{u})) = \\ (\bar{w},f) + (\tau(\boldsymbol{\alpha} \cdot \nabla \bar{w} + \kappa \nabla^2 \bar{w}),f)
\end{align}
which is precisely the stabilized formulation proposed in \cite{Franca}.  In \cite{Franca}, the stabilization parameter is defined based on error analysis considerations.  In the present work, we have derived a stabilization parameter based on a decomposition of the generalized finite element method.  To further illustrate the relationship between the generalized finite element method of advection--diffusion, consider the structure of $\tau$ in equation \eqref{Eqn:Define_Tau}, which is identical to the structure of the stabilization parameter derived using the Hughes variational multiscale method in \cite{Masud2}.  In \cite{Masud2}, for $\boldsymbol{N}'$, the authors employ bubble functions that vanish on the element boundary. The nature of the bubble functions preserves convergence as the bubble functions do not interfere with the nodal values.  In the generalized finite element framework, $\boldsymbol{N}'$ does not vanish on the element boundary, but is defined over the entire domain.  In this case, convergence is ensured based on the qualities of the partition of unity \cite{Melenk}.

	Clearly, casting the generalized finite element method into a variational multiscale framework elucidates its apparent similarities with stabilized methods.
%
\section{NUMERICAL RESULTS}
In this section, we present some numerical examples that illustrate the potential of the generalized finite element method for producing stable solutions to the advection--diffusion equation for high $Pe^h$.  Due to the sharp corners local to the enriched elements, a 100-point gauss integration scheme is used to integrate the stiffness terms in equation \eqref{Eqn:AD_Matrix_Functional_Define} for any element with enriched nodes.  Also, for the examples that enforce the boundary conditions weakly, the boundary integrals are integrated using two-point trapezoidal numerical integration.
\subsection{Linear elements in one-dimension}
\begin{figure}[htb!]
	\centering
  \subfigure[]{\includegraphics[scale=0.35]{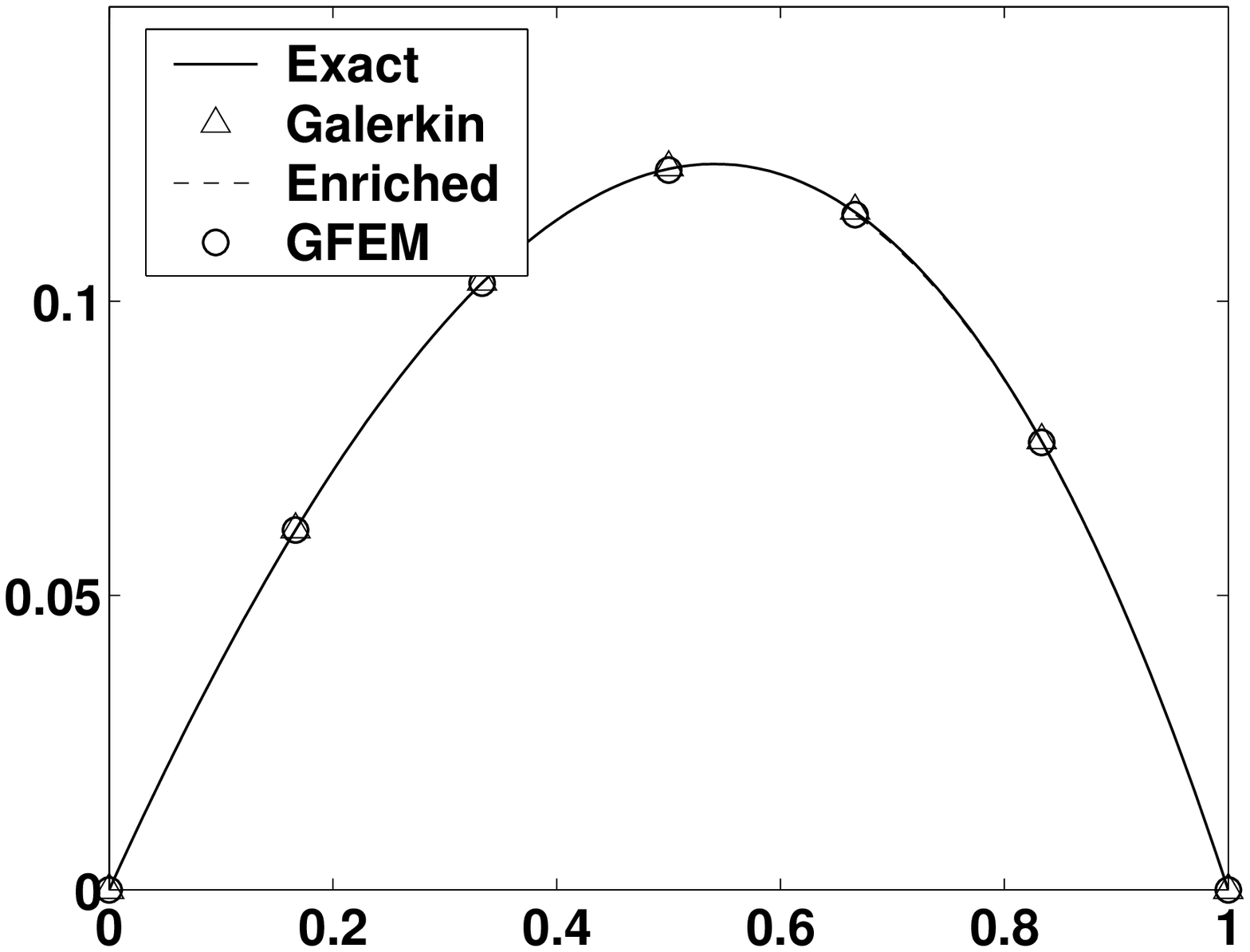}}
  \subfigure[]{\includegraphics[scale=0.35]{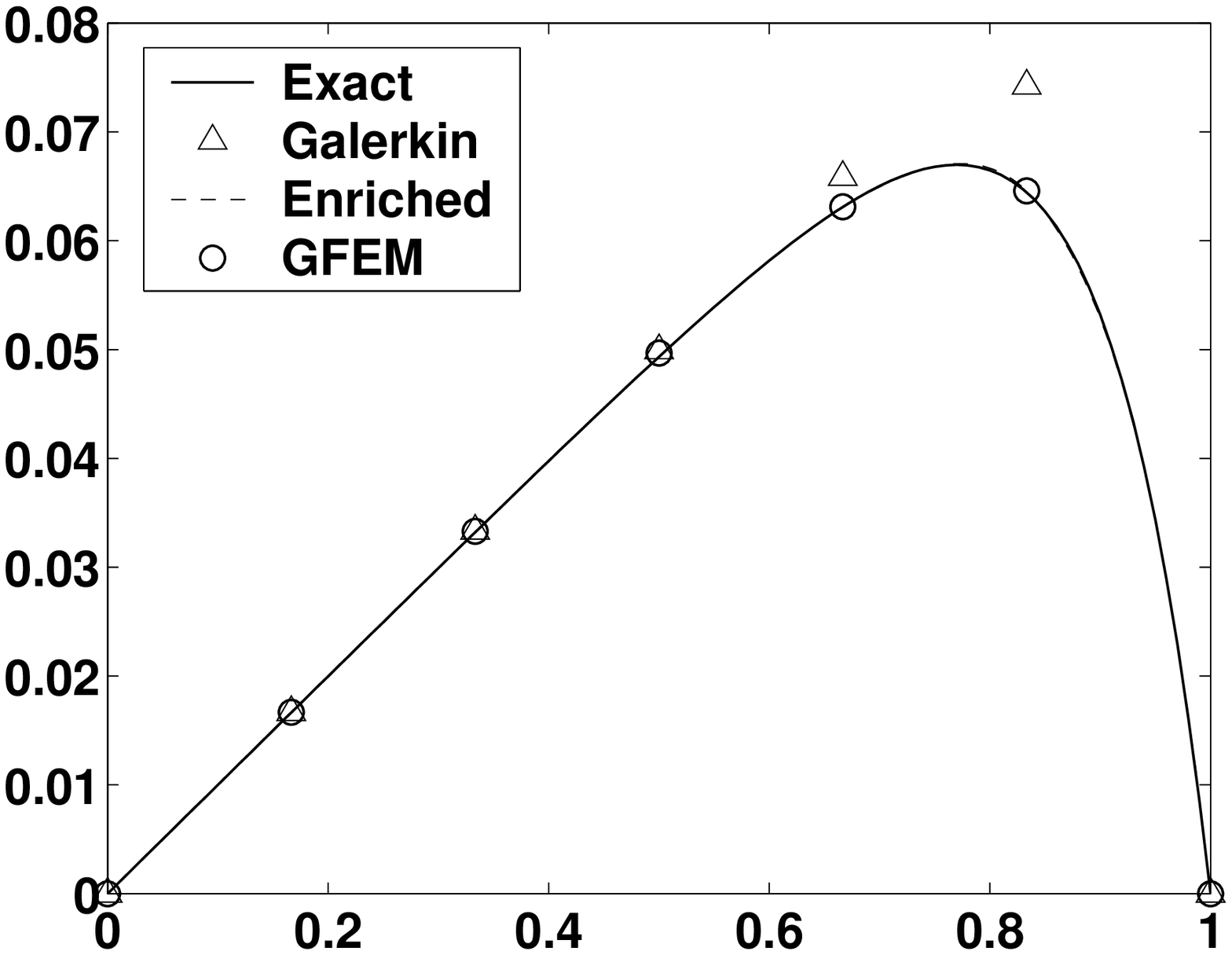}}
  \subfigure[]{\includegraphics[scale=0.35]{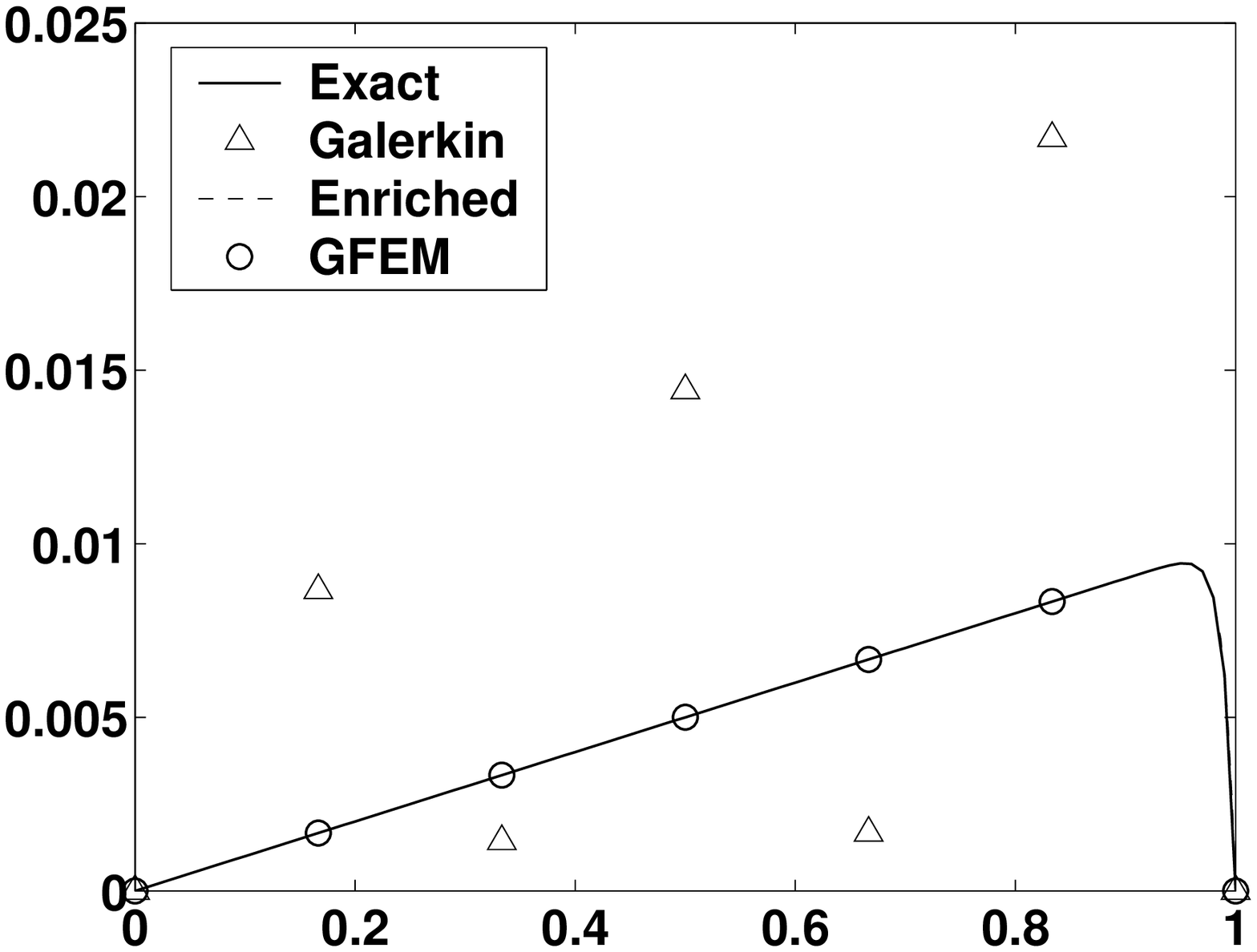}}
  \subfigure[]{\includegraphics[scale=0.35]{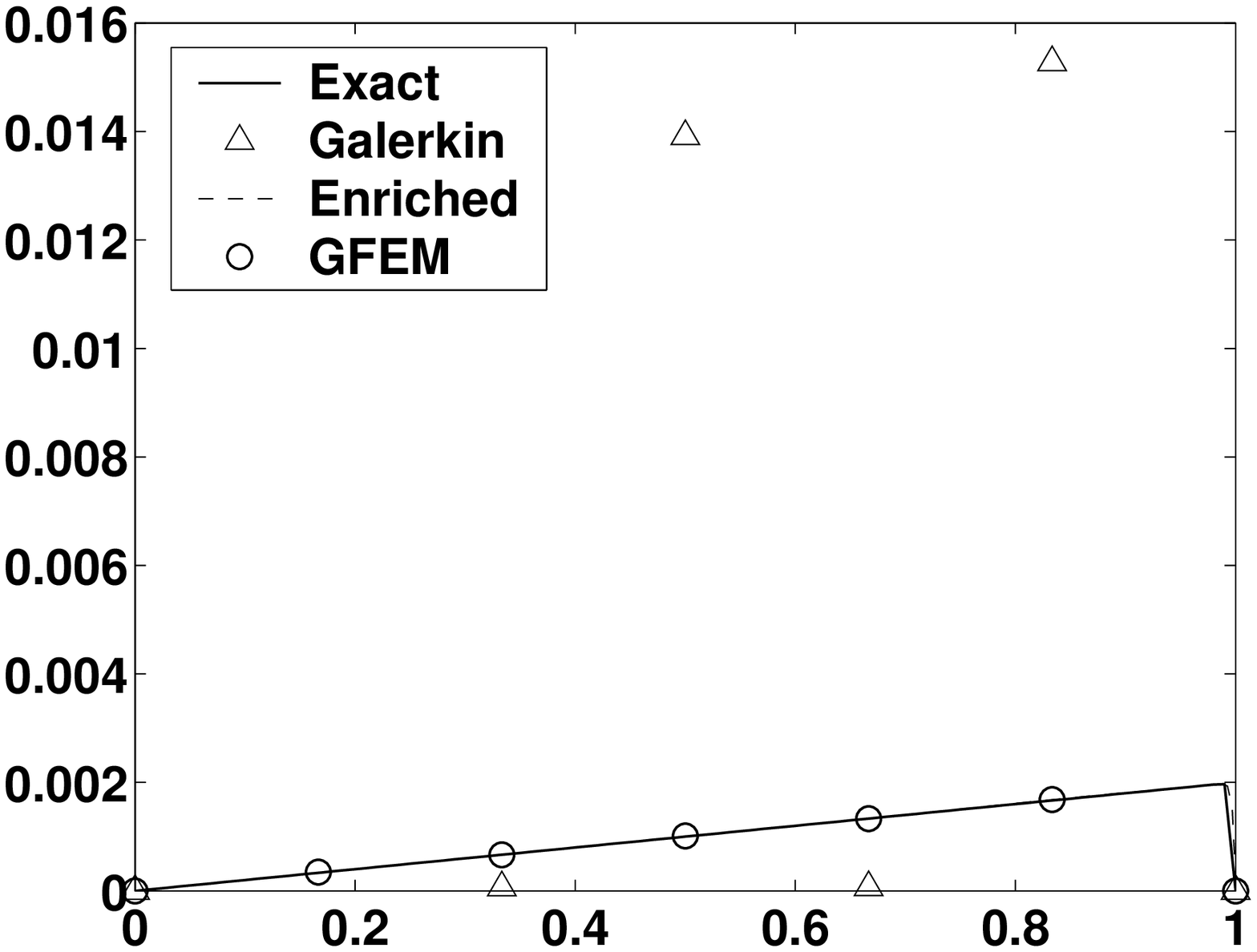}}
	\caption{Finite element solutions to the advection diffusion equation for various $Pe^h$ (a) $Pe^h = 0.16$ (b) $Pe^h = 1.66$ (c) $Pe^h = 16.66$ (d) $Pe^h = 83.33$.}
	\label{fig:1D_AD_FEM_Sols}
\end{figure}

\begin{figure}[htb!]
	\centering
  \includegraphics[scale=0.45]{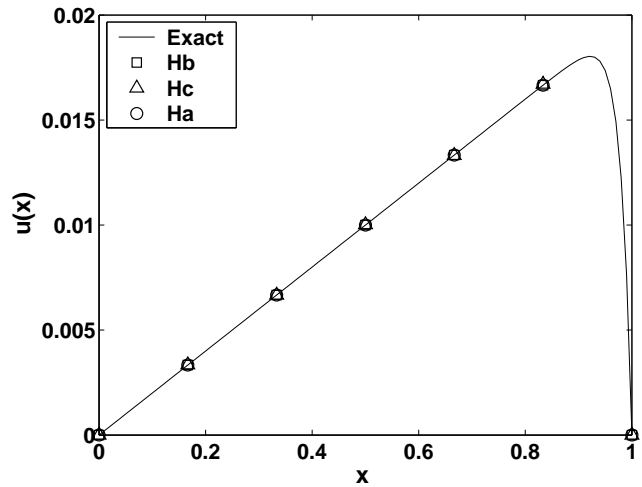}
	\caption{Generalized finite element solutions to the advection diffusion equation for $Pe^h = 5$.}
	\label{fig:1DCompare}
\end{figure}

The first example we present is the advection-diffusion equation in one-dimension, with a constant advection speed $\alpha$ and diffusivity  $\kappa$, and a unit source.
\begin{align}
   \alpha \frac{\partial u}{\partial x} - \kappa \frac{\partial^2 u}{\partial x^2} &= 1 \quad \mathrm{in} \quad \Omega \\
   u &= 0 \quad \mathrm{on} \quad \Gamma
\end{align}
The problem is solved on a unit domain discretized into six equal length elements.  The results are shown in Figure \ref{fig:1D_AD_FEM_Sols}.  Notice that the generalized finite element method performs well for high $Pe$.  Also notice that the solution between the nodes is indistinguishable from the exact solution for the generalized finite element method.  The results shown in Figure \ref{fig:1D_AD_FEM_Sols} were obtained by enriching the last two nodes with enrichment function $\mathcal{H}_b$.  To show the equivalence of enriching with  $\mathcal{H}_a$, $\mathcal{H}_b$, or $\mathcal{H}_c$ a comparison is shown in Figure \ref{fig:1DCompare}.  Since enrichment functions $\mathcal{H}_a$ and $\mathcal{H}_c$ do not automatically vanish on the boundary, the boundary conditions are enforced weakly in those cases.
\subsection{Linear elements in two-dimensions}
\begin{figure}[htb!]
	\centering
  \includegraphics[scale=0.35]{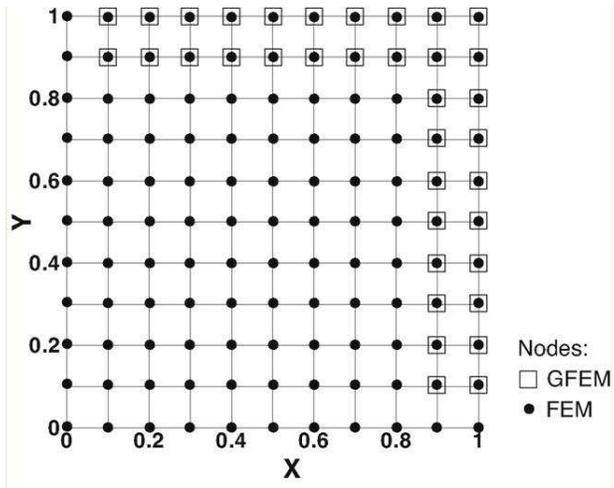}
	\caption{Finite element mesh for the 2D constant coefficients problem showing the enriched, partially enriched, and unenriched elements.}
	\label{fig:2DConstantMesh}
\end{figure}

\begin{figure}[htb!]
	\centering
  \subfigure[]{\includegraphics[scale=0.25]{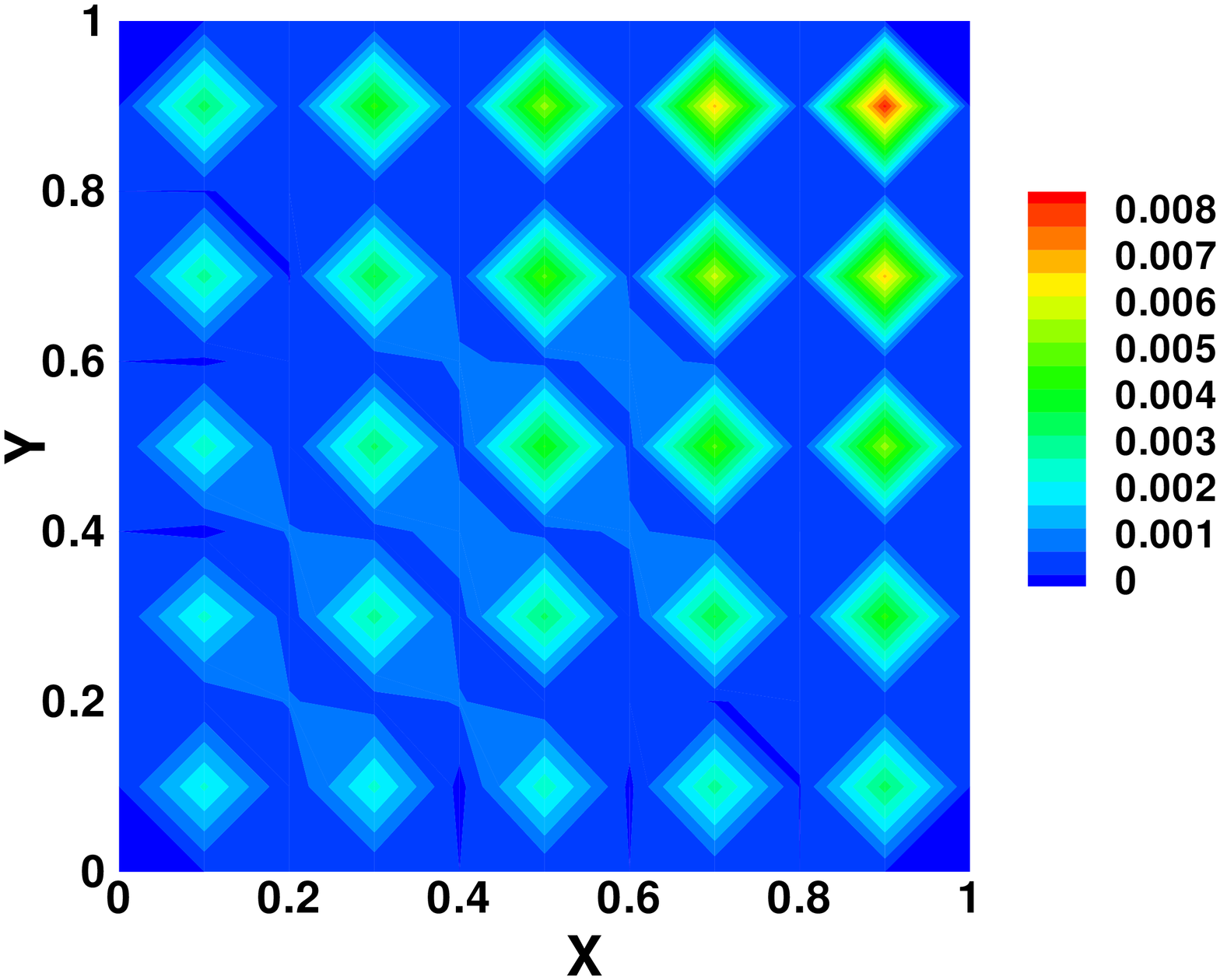}}
  \subfigure[]{\includegraphics[scale=0.25]{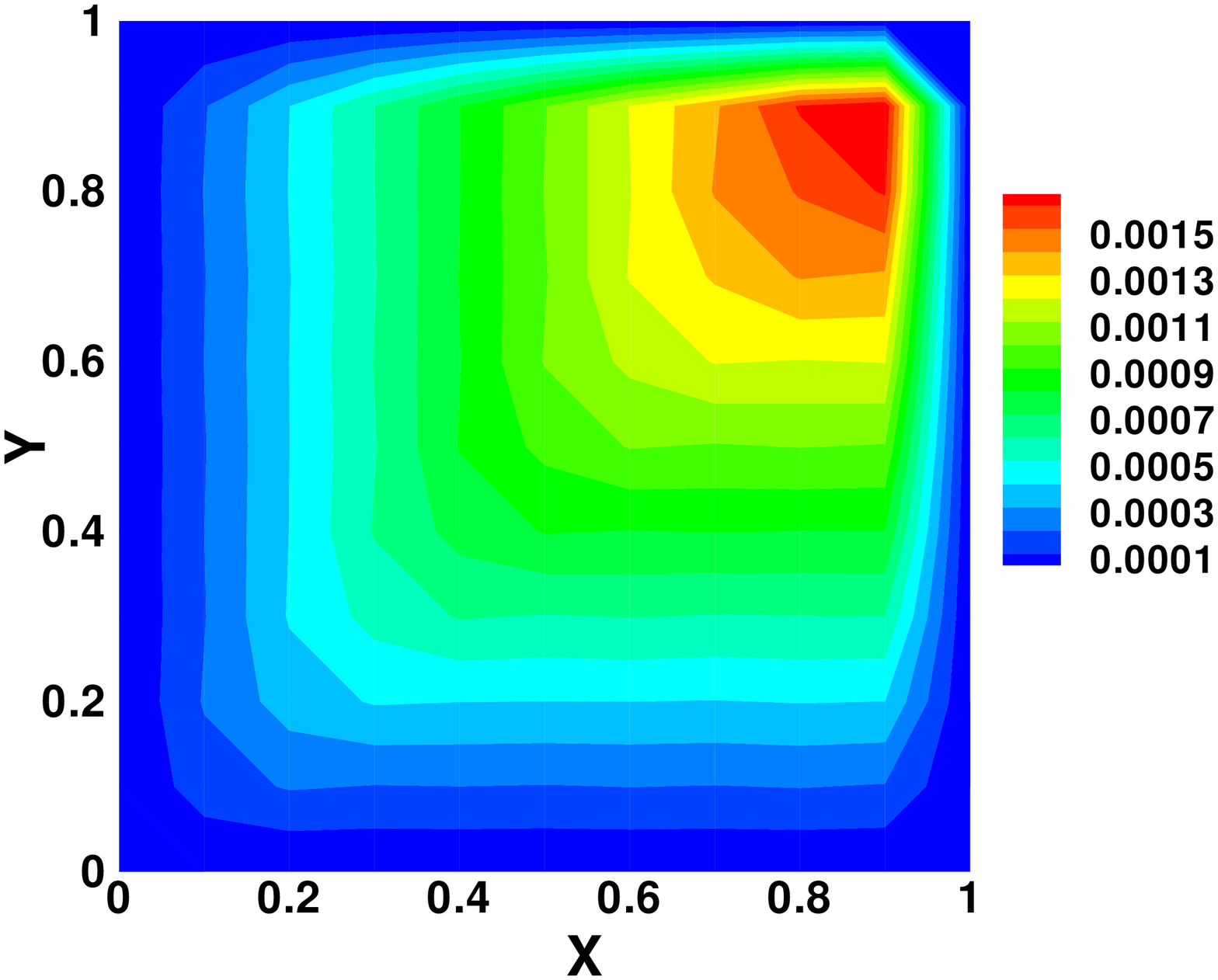}}
  \subfigure[]{\includegraphics[scale=0.25]{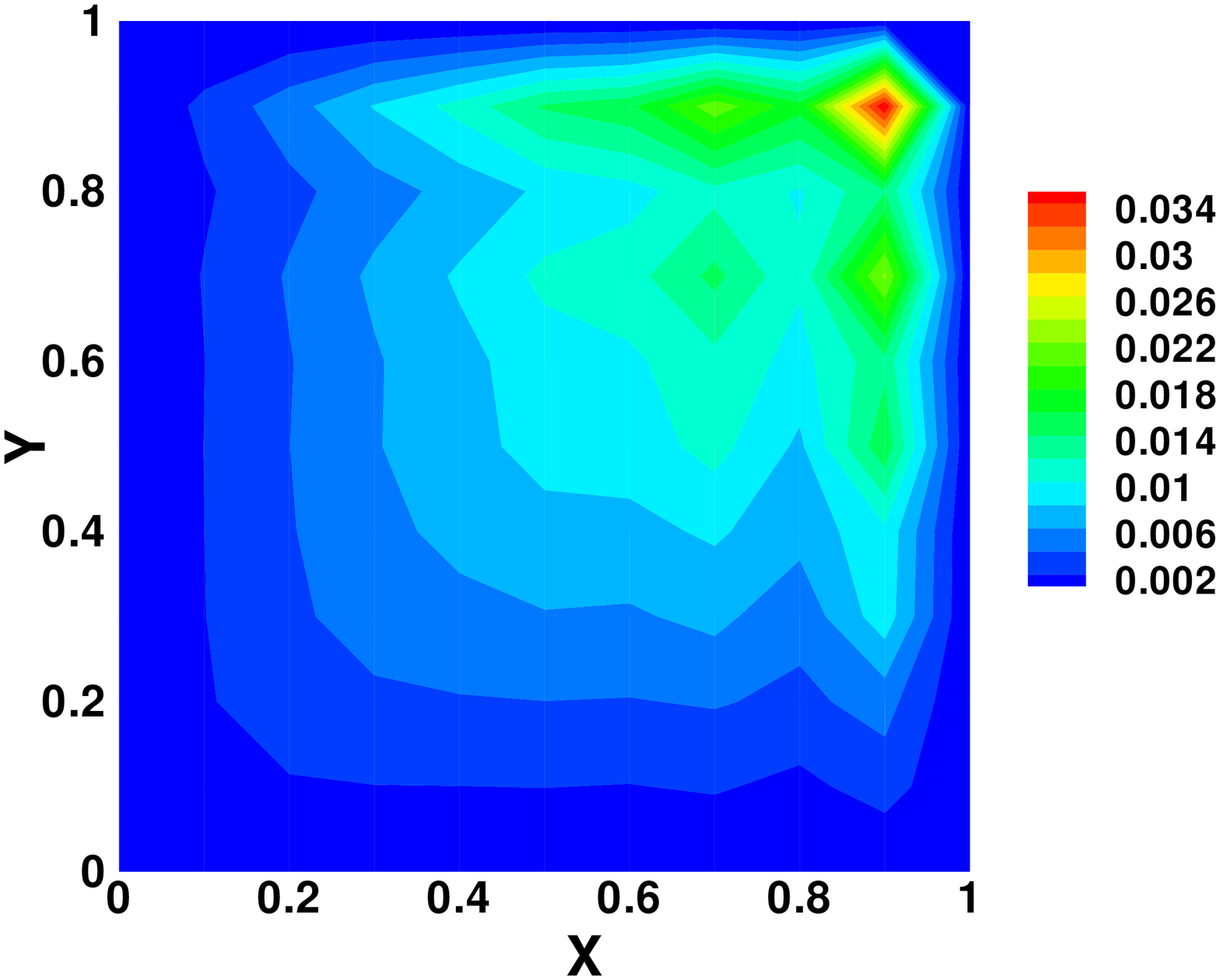}}
	\caption{Results for 2D advection diffusion equation for $Pe^h = 35$ (a) Galerkin finite element solution (b) generalized finite element solution using $\mathcal{H}_b^2$ (b) Galerkin solution using exponential finite elements.}
	\label{fig:2DConstantCoefs}
\end{figure}

The same enrichment function used in one-dimension is also applicable in two-dimensions.  To illustrate  this point, we present the results from the advection-diffusion equation applied to a unit square domain, with constant advection speed and diffusivity, and a unit source term.  The mesh used for this example is shown in Figure \ref{fig:2DConstantMesh} along with the location of the enriched nodes.  The results shown in Figure \ref{fig:2DConstantCoefs} (b) were produced by enriching with the following function
\begin{align}
	\mathcal{H}_b^2(\boldsymbol{x}) = (1 - \frac{1 - e^{\gamma x}}{1 - e^{\gamma}})(1 - \frac{1 - e^{\gamma y}}{1 - e^{\gamma}})
\end{align}
with $\gamma = \alpha/\kappa$.  Figure \ref{fig:2DConstantCoefs} (a) shows the standard Galerkin solution to the problem, which is highly oscillatory.  Figure \ref{fig:2DConstantCoefs} (c) shows a comparison with results obtained by using the standard Galerkin method with exponential shape functions rather than linear elements.
\subsection{Linear elements in one and two-dimensions using the global-local approach}
\begin{figure}[htb!]
	\centering
  \includegraphics[scale=0.45]{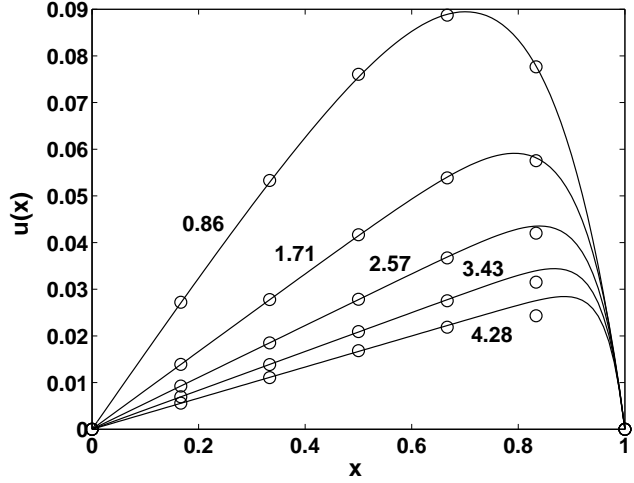}
	\caption{Generalized finite element solution to the 1D advection-diffusion equation using global-local enrichment for various $Pe^h$.}
	\label{fig:GlobalLocal1D}
\end{figure}

\begin{figure}[htb!]
	\centering
	\subfigure[]{\includegraphics[scale=0.25]{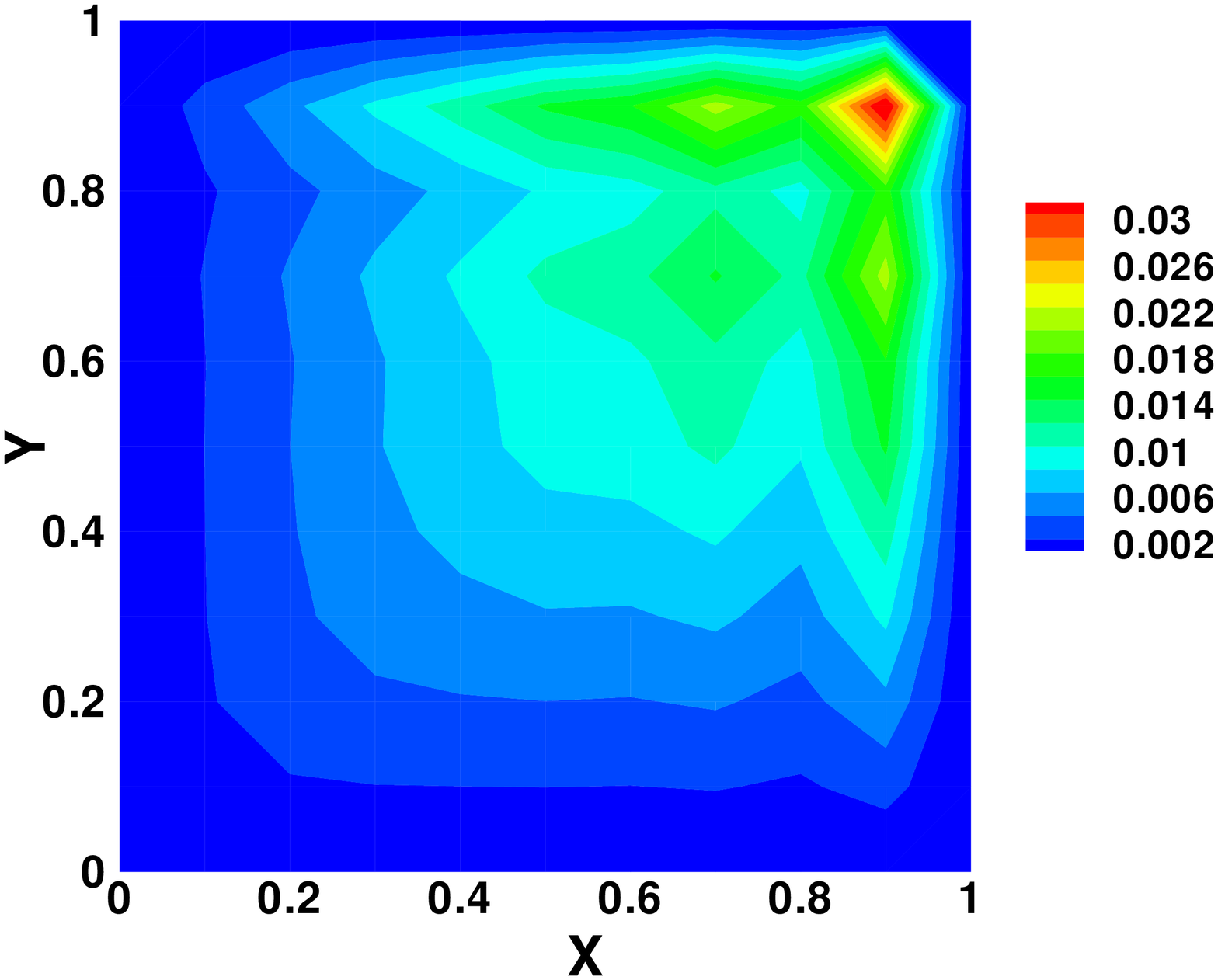}}
  \subfigure[]{\includegraphics[scale=0.25]{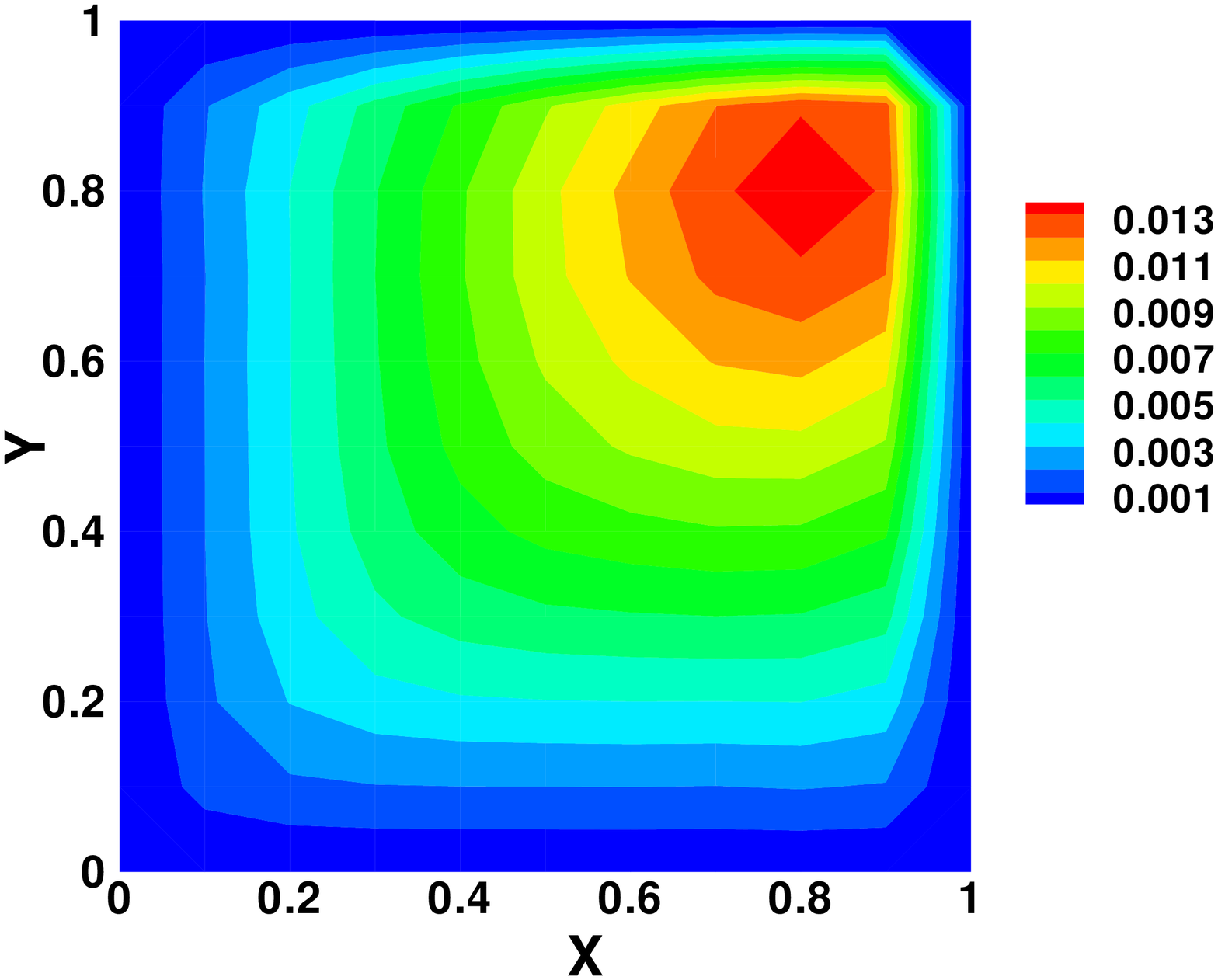}}
	\caption{2D constant coefficients on a regular domain: (a) Galerkin finite element solution (b) generalized finite element solution using global-local approach, enforcing Dirichlet boundary conditions weakly, $Pe^h = 7, \lambda = 10^6$.}
	\label{fig:LocalGlobal2D}
\end{figure}

Because the proper enrichment function is not always an obvious choice for every problem, we present some numerical examples for which the solution to the problem at low $Pe^h$ is used for the enrichment function.  Figure \ref{fig:GlobalLocal1D} shows the results from using the global-local approach to solve the same 1D problem described above with a unit source, vanishing $u$ on the boundary, on a unit domain.  Whereas using enrichment function $\mathcal{H}_b$ produces accurate results for high $Pe^h$ (as shown in Figure \ref{fig:1D_AD_FEM_Sols}), the global-local approach provides accurate solutions to $Pe^h \approx 3$.  In two-dimensions, a smooth solution is obtained for $Pe^h = 7$, shown in Figure \ref{fig:LocalGlobal2D}.
\subsection{Thermal boundary layer}
\begin{figure}[htb!]
	\centering
  \includegraphics[scale=0.6]{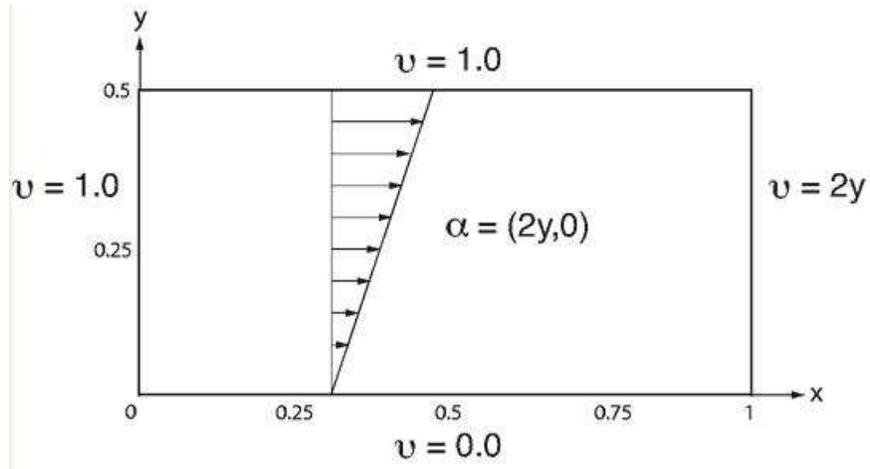}
	\caption{Thermal boundary layer problem: geometry and boundary conditions.}
	\label{fig:ThermalFigure}
\end{figure}

\begin{figure}[htb!]
	\centering
  \includegraphics[scale=0.25]{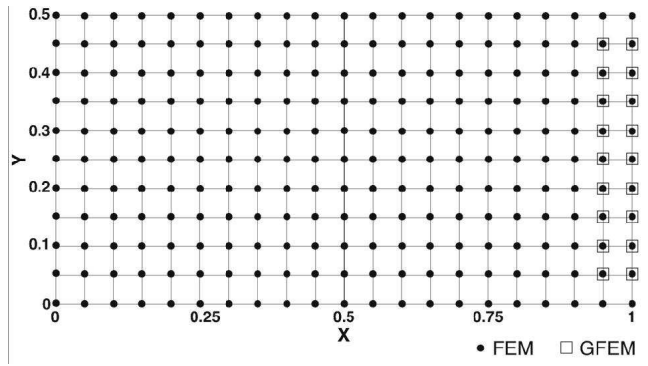}
	\caption{Finite element mesh for the thermal boundary layer problem showing the enriched, partially enriched, and unenriched elements.}
	\label{fig:ThermalMesh}
\end{figure}

\begin{figure}[htb!]
	\centering
	\subfigure[]{\includegraphics[scale=0.25]{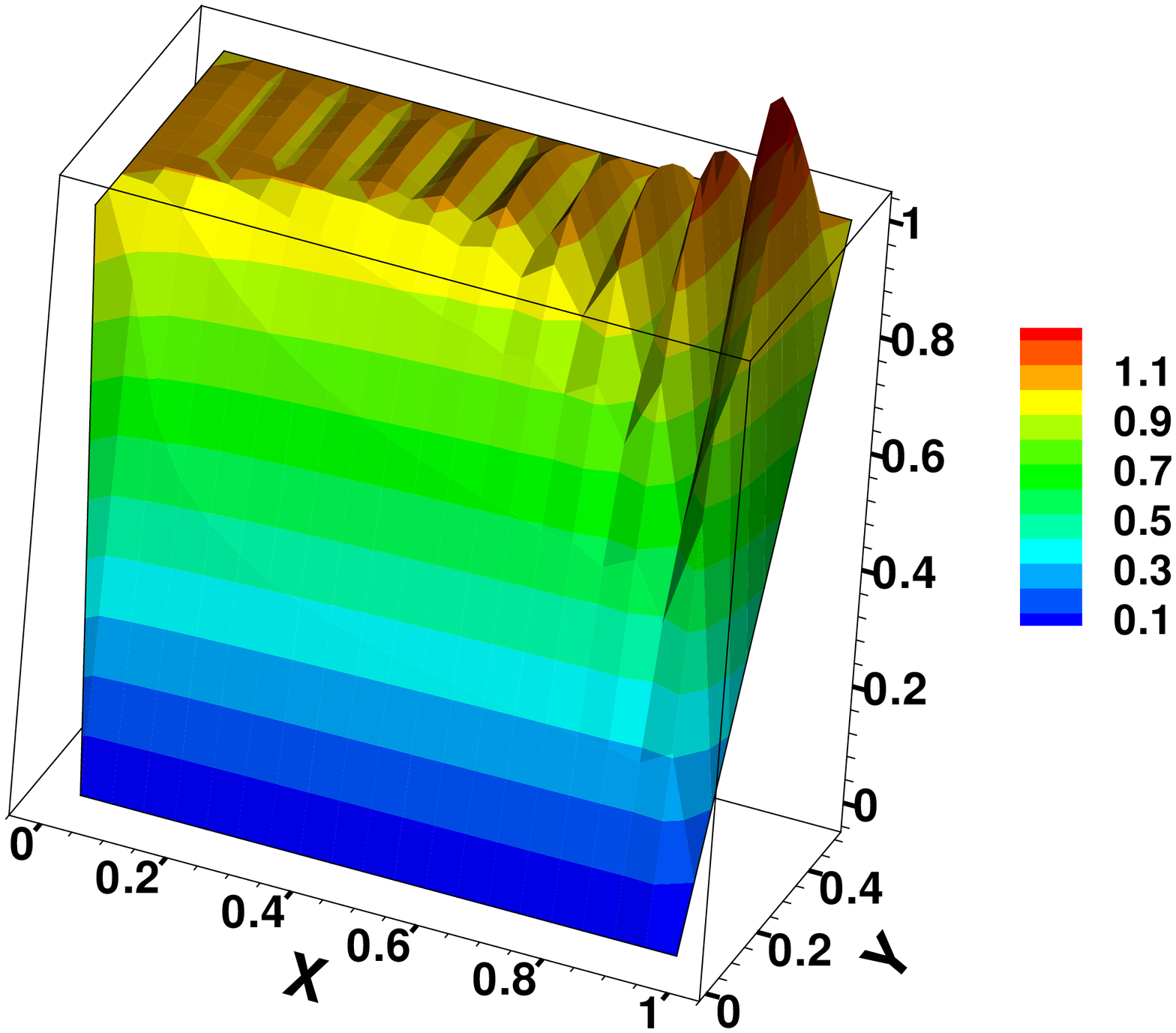}}
  \subfigure[]{\includegraphics[scale=0.25]{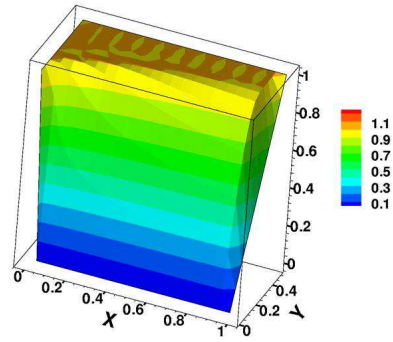}}
  \subfigure[]{\includegraphics[scale=0.25]{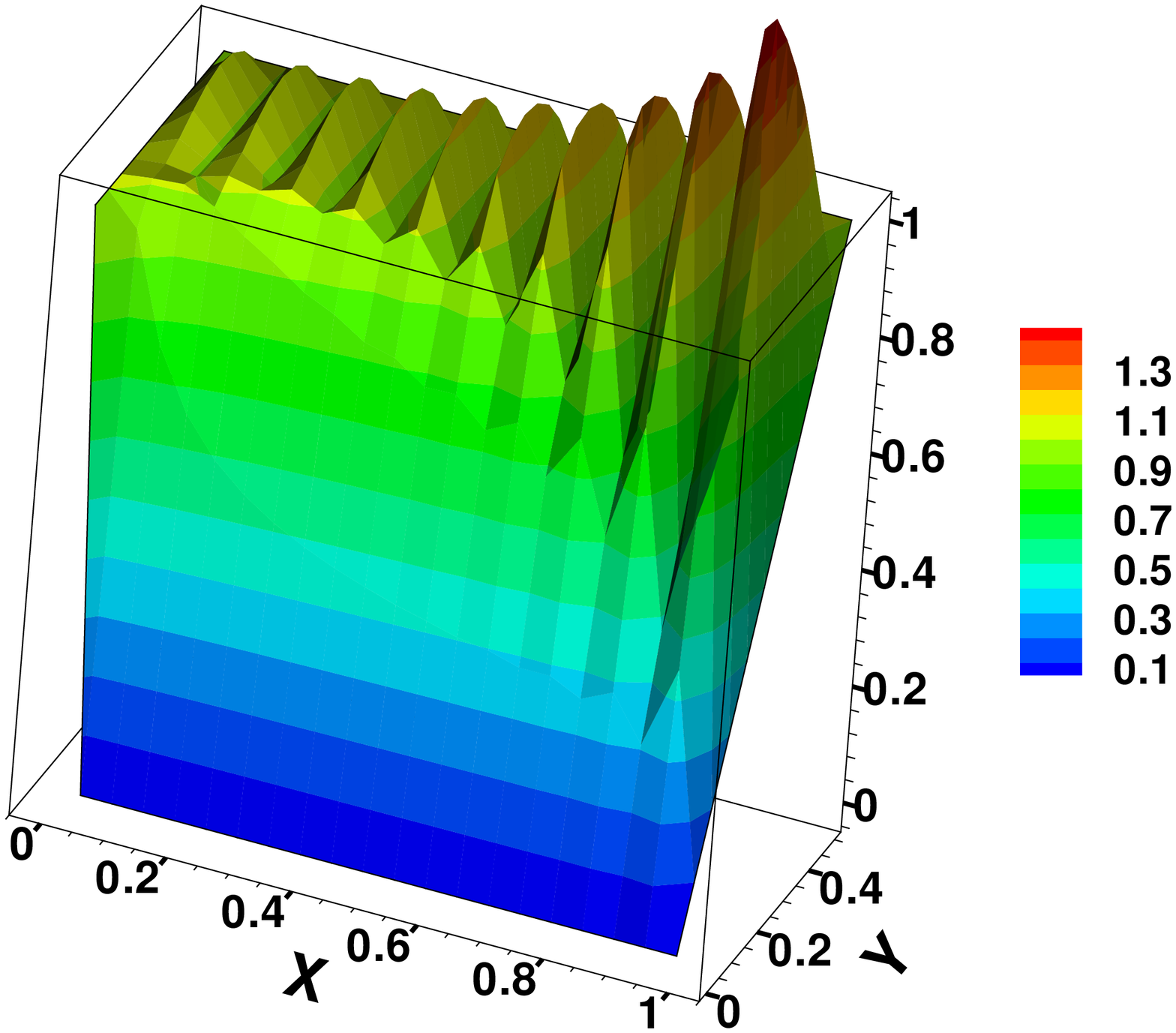}}
  \subfigure[]{\includegraphics[scale=0.25]{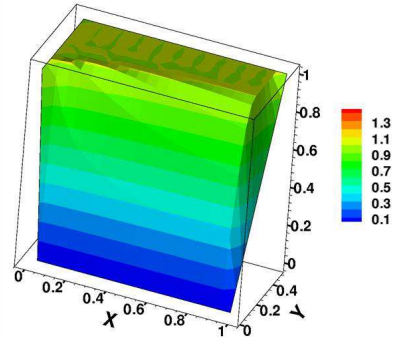}}
	\caption{Thermal boundary layer: (a) Galerkin finite element solution $Pe^h = 12.5$ (b) generalized finite element solution using $\mathcal{H}_c$, enforcing Dirichlet boundary conditions weakly with $\lambda = 10^{10}$, $Pe^h = 2.5$ (c) Galerkin finite element solution $Pe^h = 25$ (d) generalized finite element solution using $\mathcal{H}_c$, enforcing Dirichlet boundary conditions weakly with $\lambda = 10^{10}$, $Pe^h = 25$.}
	\label{fig:Plates3D}
\end{figure}

The thermal boundary layer problem was first analyzed in \cite{Franca} and represents the simulation of fully developed flow between two plates, the top plate moving with a unit velocity in the x-direction.  The domain and boundary conditions are shown in Figure \ref{fig:ThermalFigure}.  The highest global $Pe$ for this problem is $\approx 700$, whereas the element $Pe^h$ using the elements nearest the top plate is 25. The mesh with enriched nodes along the outflow boundary is shown in Figure \ref{fig:ThermalMesh}.  
	As shown in Figure \ref{fig:Plates3D}, smooth solutions are obtained using enrichment function $\mathcal{H}_c$, enforcing the boundary conditions weakly.  Strongly enforcing the boundary conditions causes spurious oscillations as predicted in \cite{Cho}.  Notice that the Galerkin solution has wiggles that propagate through the entire domain.  The results compare well with stabilized methods presented in \cite{Mulder, Franca} for the SUPG method, and \cite{Masud2} for GLS and the Hughes variational multiscale method, further illustrating the close relationship between stabilized and generalized finite element methods.
%
%

\section{CONCLUSIONS}
We have introduced a generalized finite element formulation for the advection--diffusion equation with stabilizing properties.  Since the methods is highly dependent on the enrichment function selected, we have provided various enrichment functions with different qualities and shown an equivalence between them.  We have also explored the relationship between generalized finite element methods and stabilized methods. The similarity between GLS-type methods and stabilized methods using bubble functions has been extended to the generalized finite element context by decomposing the problem in a variational multiscale manner.  The numerical examples presented show that the method is applicable to high $Pe$ flows and works for problems with non-constant boundary conditions.  The examples also illustrate the method's similarity to stabilized methods.  Finally, this work emphasizes the need for a broader interpretation of stabilized methods for flow problems including a more comprehensive theory from which stabilized and generalized methods emanate from.
%

\section*{ACKNOWLEDGMENTS}
The research reported herein was supported by the Computational Science and Engineering Fellowship (D. Z. Turner) and The Department of Energy (K. B. Nakshatrala) through a SciDAC-2 project (Grant No. DOE DE-FC02-07ER64323). This support is gratefully acknowledged. The opinions expressed in this paper are those of the authors and do not necessarily reflect that of the sponsor.

\bibliographystyle{unsrt}
\bibliography{../../DISSERTATION/References/references}
\newpage

\section*{FIGURES}
%
%
%
%
%

%
%
%
%
%

%
\end{document}